\documentclass[12pt]{article}
\usepackage{latexsym}
\usepackage{amsmath, amscd, amssymb}
\usepackage{eepic,epic}
\usepackage{srcltx}
 \usepackage{graphicx}
\usepackage{graphics,psfrag}
\begin{document}
\baselineskip 16pt

\def \Aut {\mbox{\rm Aut\,}}
\def \Cay {\mbox{\rm Cay\,}}
\newfont{\sBbb}{msbm7 scaled\magstephalf}
\newcommand{\BC}{\mbox{\Bbb C}}
\newcommand{\BR}{\mbox{\Bbb R}}
\newcommand{\BZ}{\mbox{\Bbb Z}}
\newcommand{\sZ}{\mbox{\sBbb Z}}

\newcommand{\Iso}{\mbox{Iso}\,}
\newcommand{\Fi}{\mbox{Fix}\,}
\renewcommand{\labelenumi}{\theenumi}
\newcommand{\qed}{\mbox{\raisebox{0.7ex}{\fbox{}}}}
\newtheorem{theorem}{Theorem} \rm %[section] \rm
\newtheorem{problem}{Problem}
\newtheorem{defin}{Definition}
\newtheorem{lemma}[theorem]{Lemma}
\newtheorem{prop}[theorem]{Proposition}
\newtheorem{conj}{Conjecture}
\newtheorem{op}{Open Problem}
\newtheorem{example}{Example} %[section]
\newtheorem{note}{Note}
\newtheorem{remark}{Remark}
\newtheorem{corollary}[theorem]{Corollary}
\newtheorem{construction}{Construction}
\newtheorem{proposition}{Proposition}
\newtheorem{question}{Question}
\newenvironment{pf}{\medskip\noindent{Proof:}
  \hspace{-.5cm}      \enspace}{\hfill \qed \newline \smallskip}
\newenvironment{pflike}[1]{\medskip\noindent {#1}
   \enspace}{\medskip}
\newenvironment{defn}{\begin{defin} \em}{\end{defin}}
\newenvironment{nt}{\begin{note} \em}{\end{note}}
\newenvironment{rem}{\begin{remark} \em}{\end{remark}}
\newenvironment{examp}{\begin{example} \rm}{\end{example}}
\newcommand{\lra}{\longrightarrow}
\newcommand{\vect}[2]{\mbox{$({#1}_1,\ldots,{#1}_{#2})$}}
\newcommand{\comb}[2]{\mbox{$ \left( \begin{array}{c}
        {#1} \\ {#2} \end{array}\right)$}}
\newcommand{\tve}[1]{\mbox{$\mathbf{\tilde {#1}}$}}
\newcommand{\bve}[1]{\mbox{$\mathbf{{#1}}$}}
\setlength{\unitlength}{12pt}
\renewcommand{\labelenumi}{(\theenumi)}
\def\mod{\hbox{\rm mod }}

\title{Chromatic-choosability of the power of graphs}
 \author{Seog-Jin Kim\thanks{Department of   Mathematics Education,
Konkuk University, Korea,   Email: skim12@konkuk.ac.kr.
  } \and
 Young Soo Kwon \thanks{Department of Mathematics, Yeungnam University, Korea,    Email: ysookwon@ynu.ac.kr }
  \and
 Boram Park\thanks
 {National Institute for Mathematical Sciences, Daejeon 305-811, Korea. E-mail: borampark@nims.re.kr}}

\date{}
\maketitle

\begin{abstract}
The $k$th power $G^k$ of a graph $G$ is the graph defined on $V(G)$ such that two vertices $u$ and $v$ are adjacent in $G^k$ if the distance between $u$ and $v$ in $G$ is at most $k$.  Let $\chi(H)$ and $\chi_l(H)$ be the chromatic number and the list chromatic number of $H$, respectively. A graph $H$ is called {\em chromatic-choosable} if $\chi_l (H) = \chi(H)$.
It is an interesting problem to find graphs that are chromatic-choosable.
A natural question raised by Xuding Zhu \cite{Zhu2012} is whether there exists a constant integer $k$ such that $G^k$ is chromatic-choosable for every graph $G$.

Motivated by the List Total Coloring Conjecture, Kostochka and Woodall \cite{KW2001} asked whether $G^2$ is chromatic-choosable for every graph $G$.
Kim and Park \cite{KP2013} answered the Kostochka and Woodall's question in the negative by
finding a family of graphs whose squares are complete multipartite graphs with partite sets of equal and unbounded size.
In this paper, we answer Zhu's question by showing that
for every integer $k \geq 2$,  there exists a graph $G$ such that $G^k$ is  not chromatic-choosable. Moreover, for any fixed $k$ we show that the value $\chi_l(G^k) - \chi(G^k)$ can be arbitrarily large.

 \vskip10pt

\end{abstract}

\section{Introduction}

For any graph $G$ and for any positive integer $k$,
the $k$th power $G^k$ of a graph $G$ is the graph defined on $V(G)$ such that two vertices $u$ and $v$ are adjacent in $G^k$ if the distance between $u$ and $v$ in $G$ is at most $k$.  In particular, $G^2$ is called the \emph{square} of $G$.

A proper $k$-coloring $\phi: V(G) \rightarrow \{1, 2, \ldots, k \}$ of a graph $G$ is an assignment of colors to the vertices of $G$ so that any two adjacent vertices  receive distinct colors.
The {\em chromatic number} $\chi(G)$ of a graph $G$ is the least $k$ such that there exists a proper $k$-coloring of $G$.
A list assignment $L$ is an assignment of lists of colors to vertices.
A graph $G$ is said to be {\em $k$-choosable} if for any list $L(v)$ of size at least $k$, there exists a proper coloring $\phi$ such that $\phi(v) \in L(v)$ for every $v \in V(G)$.
The least $k$ such that $G$ is $k$-choosable is called the {\it list chromatic number} $\chi_l(G)$ of a graph $G$. Clearly $\chi_l(G) \geq \chi(G)$ for every graph $G$.

A graph $G$ is called {\em chromatic-choosable} if $\chi_l (G) = \chi(G)$. It is an interesting problem to find
 graphs that are chromatic-choosable.  There are several famous conjectures that some  classes of graphs
are chromatic-choosable including the List Coloring Conjecture  \cite{BH1985} and the List Total Coloring Conjecture \cite{BKW1997},
which say that the line graph and the total graph of any graph are chromatic-choosable, respectively.
Motivated by the List Total Coloring Conjecture, Kostochka and Woodall \cite{KW2001} proposed
the List Square Coloring Conjecture which states that $G^2$ is chromatic-choosable for every graph $G$.
It was noted in \cite{KW2001} that the List Total Coloring Conjecture is true if the List Square Coloring Conjecture is true.

If the diameter of a graph $G$ is $m$, then $G^m$ is chromatic-choosable since $G^m$ is a complete graph.
Thus given a graph $G$, there exists an integer $k$ such that $G^k$ is chromatic-choosable.
A natural question raised by Xuding Zhu \cite{Zhu2012} is whether there exists a constant $k$ such that $G^k$ is chromatic-choosable for every graph $G$.
The question looks trivially true if we do not think carefully.
But, it is not easy to answer the question whether such constant $k$ exists or not.
On the other hand, it is an interesting problem to find a class $\cal{F}$ of graphs such that for every integer $k \geq k_0$ for some constant $k_0$, $H^k$ is chromatic-chooslable for every graph $H$ in $\cal{F}$.
In this direction, Wang and Zhu \cite{WZ2010} show that $C_n^k$ is chromatic-choosable for any positive integer $k$ and for any cycle $C_n$.

The List Square Coloring Conjecture proposed the smallest $k$ such that $G^k$ is chromatic-choosable for every graph $G$ is 2.  Also it proposed that $G^{2t}$ is chromatic-choosable for any positive integer $t$.
However, recently, Kim and Park \cite{KP2013} disproved the List Square Coloring Conjecture by showing that for any prime $n$, there exists a graph $G$ such that  $G^2$ is $K_{n \star (2n-1)}$, where
$K_{n \star (2n-1)}$ denotes the complete multipartite graph with $(2n-1)$ partite sets in which each partite set has size $n$.  Note that the gap between $\chi_l (G^2)$  and $\chi(G^2)$ can be arbitrarily large since the gap between the list chromatic number and chromatic number of a complete multipartite graph can be arbitrarily large.

Since there exists a graph $G$ such that $G^2$ is not chromatic-choosable, next direction is to find the smallest $k$ such that $G^{k}$ is chromatic-choosable for every graph $G$, if such $k$ exists. In this paper, we answer Zhu's question by showing that there is no constant $k$ such that $G^k$ is chromatic-choosable for every graph $G$.  We show the following main theorem.

\begin{theorem}\label{main}
For any positive integer $k \geq 2$ and
for any positive integer $s$, there exists a graph $G$ such that
\[\chi_l(G^k) \ge \frac{10}{9} \cdot 3^{3sk-1} - 1 \ge 3^{3sk-1} = \chi(G^k).\]
\end{theorem}

This implies that for any integer $k \geq 2$, there exists a graph $G$ such that $G^k$ is not chromatic-choosable.  Moreover, since
\[\chi_l(G^k) -  \chi(G^k) \ge  (\frac{10}{9}-1)  \cdot 3^{3sk-1}  -1 =  3^{3sk-3} -1,\]
for any fixed $k$, the value $\chi_l(G^k) - \chi(G^k)$ can be arbitrarily large as the integer $s$ goes to infinity.

%%%%%%%%%%%%%%%%%%%%%%%%%%%%%%%%%%%%%%%%%%%%%%%%%%%%%%%%

\section{Construction}

In this section, first we define a Cayley graph $G_{3n}$, and study the properties of $G_{3n}$.
And then we will define a graph $H_{3n}$ which will be used in the proof of
%has the desired property as in
Theorem~\ref{main}.
First, we define the Cayley graph  $G_m$ where $m$ is an integer at least 2.

\begin{construction} \label{construction} \rm

For any positive integer $m \geq 2$, let $\mathbb{Z}_{3}^{m}$ be an (usual) additive abelian group of order $3^m$. Namely $\mathbb{Z}_{3}^{m} = \{(a_1 , a_2 , \ldots, a_{m} )  :   a_i \in \{0,1,2 \} \mbox{ for all } 1 \leq i \leq m  \}$.
Note that the identity element of $\mathbb{Z}_{3}^{m}$, denoted by ${\bf 0}$, is $(0,\ldots,0)$.

For any $i,j \in \{ 1, \ldots,m \}$ with $i \neq j$, let $x_{i,j}$ be the vector such that  the $i$th coordinate of $x_{i,j}$ is $1$, the $j$th coordinate of $x_{i,j}$ is $2$, and all other coordinates of $x_{i,j}$ are $0$.
For example, $x_{1,3} = (1,0,2,0,\ldots,0)$.
For any positive integer $m \geq 2$, let $X_m$ be the subset of $\mathbb{Z}_{3}^{m}$ such that
\begin{eqnarray*}
X_m &=& \{ x_{i,j} : 1 \leq i,j \leq m,  \mbox{ and }   i \neq j  \} .
\end{eqnarray*}
Let
\[\Gamma_m =\{ z\in \mathbb{Z}_3^m : \sum_{i=1}^{m}z_i=0\pmod{3} \},\]
where $z_i$ is the $i$th coordinate of $z$ for any $i\in \{1, \ldots,m\}$.
Now $\Gamma_m$ is a subgroup of $\mathbb{Z}_3^m$ of index $3$, since $y-z\in \Gamma_m$ for any two $y,z\in \Gamma_m$. We define a graph $G_m$ as follows:
\begin{eqnarray*}
V(G_m)&=& \Gamma_m\\
E(G_m)&=&\{ yz : y-z\in X_m\}.
\end{eqnarray*}
Note that $|V(G_m)|=|\Gamma_m|=3^{m-1}$ and
$G_m$ is a connected graph.
Actually, the graph $G_m$ is the Cayley graph on $\Gamma_m$ with the symmetric generating
set $X_m$.

\end{construction}

%%%%%%%%%%%%%%%%%%%%%%%%%%%%%%%%%%%%%%%%%%%%%%%%%%%%%%%%%%%%%%%%%%%%%%%%%%%%%%%%%%%%%%%%%%%%%%%

From now on, $G_{3n}$ denotes the graph defined in Construction \ref{construction} for $m = 3n$.
Let ${\bf a}_{3n}$, ${\bf b}_{3n}$ be %, ${\bf c}_{3n+1}$, ${\bf d}_{3n+1}$ be
 the vectors in $\mathbb{Z}_3^{3n}$ defined by
\begin{eqnarray*}
{\bf a}_{3n}&=&(1,1,\ldots,1) \\
{\bf b}_{3n}&=&(2,2,\ldots,2).
\end{eqnarray*}
First, we study basic properties of the graph $G_{3n}$ in Lemma~\ref{lem3} and Lemma~\ref{lemma_G}.

\begin{lemma}\label{lem3}
Let $n\ge 2$ be an integer.
For any vector $y \in \Gamma_{3n}$,
\[d_{G_{3n}}({\bf 0},y) \le  \frac{2\times (\text{the number of nonzero coordinates of }y)}{3}, \]
and the equality holds if and only if all nonzero coordinates of $y$ are identical.
\end{lemma}

\begin{pf}
We will prove the lemma by the induction on the number of nonzero coordinates of $y$.
If $y \neq {\bf 0}$, then $y$ has at least two nonzero coordinates.

As the basis step, we consider the cases that $y$ has two or three nonzero coordinates.
If $y$ has two  nonzero coordinates, then one can check $y=x_{i,j}$ for some $i,j \in \{1,\ldots,3n\}$, and hence the lemma holds.
If  $y$ has three nonzero coordinates, then all nonzero coordinates of $y$ are identical and $d_{G_{3n}}({\bf 0},y) = 2$. So the lemma holds.
As the induction step, suppose that the lemma holds for any $y$ which has at most $m$ nonzero coordinates ($m\ge 3$).
Now we take a vertex $y$ of $G_{3n}$  which has $(m+1)$ nonzero coordinates.
For convenience, we denote by $y_i$ the $i$th coordinate of $y$ for any $i\in \{1, \ldots,3n\}$.

\medskip

\noindent {\bf Case 1.} There exist two distinct integers $i_1, i_2  \in \{1, \ldots,3n\}$ such that $y_{i_1}=1$ and $y_{i_2}=2$

Let $w$ the vector such that $y=w+x_{{i_1,i_2}}$.  Note that $w$ has  $(m+1)-2$ nonzero coordinates.
Then by the induction hypothesis,
\[d_{G_{3n}}({\bf 0},y)\leq d_{G_{3n}}({\bf 0},w)+1  \leq \frac{2((m +1) -2)}{3} + 1 < \frac{2(m+1)}{3}.\]

\medskip

\noindent {\bf Case 2.}  All nonzero coordinates of $y$ are identical.

Since $y$ has $(m+1)$ nonzero coordinates and  $m+1\ge 4$,
there are three distinct integers $i_1, i_2, i_3 \in \{1, \ldots,3n\}$ such that $y_{i_1}=y_{i_2}=y_{i_3}\neq 0$.
First, we suppose that  $y_{i_1}=y_{i_2}=y_{i_3}=1$.
Let $w$ the vector such that  $y=w+x_{{i_1,i_2}}+x_{{i_3,i_2}}$.
Note that $w$ has $(m+1)-3$ nonzero coordinates and all of the nonzero coordinates of $w$ are identical.
Thus by the induction hypothesis,
\[d_{G_{3n}}({\bf 0},y) \leq d_{G_{3n}}({\bf 0},w)+ 2  = \frac{2(m+1-3)}{3} +2=\frac{2(m+1)}{3}.\]
Therefore, $d_{G_{3n}}({\bf 0},y)\le \frac{2(m+1)}{3}$.

Note that since the sum of all coordinates of $y$ is 0 modulo $3$ and any nonzero coordinate of $y$ is 1, the number of nonzero coordinates of $y$ is multiple of $3$. Therefore $\frac{2(m+1)}{3}$ is an integer.
Now we will show that $d_{G_{3n}}({\bf 0},y) \ge \frac{2(m+1)}{3}$ to conclude that $d_{G_{3n}}({\bf 0},y) = \frac{2(m+1)}{3}$.
%if all nonzero coordinates of $y$ are identical.

Suppose that $d_{G_{3n}}({\bf 0},y) < \frac{2(m+1)}{3}$.
Then there is a  subset $A$ of $X_{3n}$ such that
$y=\sum_{x_{i,j}\in A}x_{i,j} \text{ and }  |A| =d_{G_{3n}}({\bf 0},y)$.
We define a graph $H$ such that $V(H)=\{1,2,\ldots,3n\}$ and $E(H)=\{ ij \mid x_{i,j}\in A\}$.
Note that $|E(H)| < \frac{2(m+1)}{3}$ since $|A|< \frac{2(m+1)}{3}$.
Let $W_1$, $W_2$, $\ldots$, $W_s$ be nontrivial connected components of $H$.
Let $H_0$ be the union of  $W_1$, $W_2$, $\ldots$, $W_s$.
Since each nonzero coordinate $i$ of $y$ cannot be an isolated vertex in $H$, $i$ must belong to $V(H_0)$, and so  $|V(H_{0})| \ge  m+1$. Note that  $|E(H_0)|=|E(H)|$.
 %Next, we will show that $|V(H_{0})| \leq  m+1$  to conclude that $|V(H_{0})| =  m+1$.
%Suppose that there exists $j \in \{1,\ldots,3n \}$ such that $y_j=0$ and $j$ belongs to $V(H_0)$. Then there %are distinct integers $i_1 , i_2, i_3 \in \{1,\ldots,3n\}$ such that $x_{j, i_1 }, x_{i_2 , j } \in A$,
%$x_{j, i_1 }, x_{j , i_2 }, x_{j , i_3 } \in A$, or $x_{ i_1, j }, x_{i_2, j }, x_{i_3 , j} \in A$.
%Since
%\[ x_{j, i_1 } +  x_{i_2 , j } = x_{i_2 , i_1},~~x_{j, i_1 } + x_{j , i_2 } +  x_{j , i_3 } = x_{i_1 , i_2 } +  %x_{i_1 , i_3 },~~\mbox{and}~~x_{ i_1, j }+ x_{i_2, j }+  x_{i_3 , j}=x_{ i_1, i_3 }+ x_{i_2, i_3 },\]
%there is  a  subset $A'$ of $X_{3n}$ such that   $y=\sum_{x_{i,j}\in A'}x_{i,j}$ and $|A'| < |A|$. This %implies $d_{G_{3n}}({\bf 0},y) \le |A' | < |A|$, which is a contradiction to (\ref{definition-A}).
%Thus if $y_j =0$ for some $j \in \{1,\ldots,3n \}$, then $j \notin V(H_0)$.
%Therefore $|V(H_{0})| \leq  m+1$, and  we conclude that $|V(H_{0})| =  m+1$.
On the other hand, it is true that  $|E(H_0)|+s \ge |V(H_0)|$.
Thus
\[ s \ge |V(H_0)| -|E(H)|>(m+1)-\frac{2(m+1)}{3}=\frac{m+1}{3}.\]

Next, we will show that $|E(W_i)| \ge 2$ for each $i \in \{1, \ldots,s\}$.
It is clear $|E(W_i)| \ge 1$ since $W_i$ is nontrivial.
If  $|E(W_i)|= 1$ for some $i \in \{1,\ldots,s\}$, say $E(W_i)=\{jk\}$, then
one of $y_j$ and $y_k$ must be 2, a contradiction to the assumption that any nonzero coordinate of $y$ is 1.
Thus we conclude that  $|E(W_i)| \ge 2$ for all $i \in \{1,\ldots,s\}$.
Therefore we have
\[ |A| =|E(H)|=|E(H_0)|=\sum_{i=1}^{s} |E(W_i)| \ge 2s >2\cdot \frac{m+1}{3},\]
which is a contradiction. Thus $d_{G_{3n}}({\bf 0},y) \ge \frac{2(m+1)}{3}$, and hence $d_{G_{3n}}({\bf 0},y) = \frac{2(m+1)}{3}$.
Therefore the lemma holds for the vector $y$.
For the case that $y_{i_1}=y_{i_2}=y_{i_3}=2$, one can show the lemma holds by a similar argument.

Therefore by Case 1 and Case 2, the lemma holds for the vector $y$.
\end{pf}

Define a relation $\sim$ on $V({G}_{3n})$ by
\[ y\sim z \quad \text{ if and only if }\quad y-z \in\{{\bf 0}, {\bf a}_{3n},{\bf b}_{3n}\}.\]
Then one can check the relation $\sim$ is an equivalence relation. For each vertex $y\in V(G_{3n})$, the equivalent class $[y]$ containing $y$ has three elements,
\[[y]=\{ y, y+{\bf a}_{3n}, y+{\bf b}_{3n}\}.\]

Note that $K_{n \star r}$ denotes the complete multipartite graph with $r$ partite sets in which each partite set has size $n$.

\begin{lemma}\label{lemma_G}
For any integer $n\ge 2$,
$G_{3n}^{2n-1}$ is the complete multipartite graph $K_{3\star 3^{3n-2}}$ in which partite sets are the equivalent classes obtained by the relation $\sim$.
\end{lemma}
\begin{pf}
By Lemma~\ref{lem3}, for any two vertices $y$ and $z$ in $G_{3n}$,
  $d_{G_{3n}}(y,z)\le 2n$ and the equality holds only when
$y-z \in\{ {\bf a}_{3n},{\bf b}_{3n}\}$. Therefore the lemma holds.
\end{pf}

It is proved in \cite{K2000} that the complete multipartite graph $K_{3\star 3^{3n-2}}$ is not chromatic-choosable.  Thus
Lemma~\ref{lemma_G} implies that the $(2n-1)$th power of $G_{3n}$ is not chromatic-choosable.
This implies that
for every odd integer $2k-1$, there exists a graph $G$ such the $(2k-1)$th power of $G$, denoted by $G^{2k-1}$,  is not chromatic-choosable.

%%%%%%%%%%%%%%%%%%%%%%%%%%%%%%%%%%%%%%%%%%%%%%%%%%%%%%%%%%%%%%%%%%%%%%%%%%%%%%%%%%%%%%%
Next, we define another graph $H_{3n}$ which will be used in the proof of Theorem \ref{main}. The \textit{cartesian product} of $G$ and $H$, denoted by  $G\square H$, is the graph with the vertex set $V(G)\times V(H)$ such that for two vertices $(g_1,h_1),(g_2,h_2)\in V(G)\times V(H)$,
 $(g_1,h_1)$ and $(g_2,h_2)$ are adjacent in $G\square H$ if and only if
either $ h_1=h_2$ and  $g_1g_2\in E(G)$, or $g_1=g_2$ and $h_1h_2\in E(H)$.

\begin{construction} \label{construction-two}
For any positive integer $n$, let $H_{3n}$ be the cartesian product $G_{3n} \square K_3$ of $G_{3n}$ and $K_3$.
%See Figure~\ref{fig2} for an illustration.
\end{construction}

%%%%%%%%%%%%%%%%%%%%%%%%%%%%%%%%%%%%%%%%%%%%%%%%%%%%%%%%%%%%%%%%%%%%%%%%%%%%%%%%%%%%
\section{Proof of Theorem~\ref{main}}

For two graphs $G$ and $H$, the \textit{lexicographic product} of $G$ and $H$, denoted by $G[H]$, is the graph with the vertex set $V(G)\times V(H)$ such that
for two vertices $(g_1,h_1),(g_2,h_2)\in V(G)\times V(H)$, $(g_1,h_1)$ and $(g_2,h_2)$ are adjacent in $G[H]$ if and only if
either $g_1g_2\in E(G)$, or $g_1=g_2$ and $h_1h_2\in E(H)$.

We will show that the $(2n)$th power of $H_{3n}$, denoted by $H_{3n}^{2n}$,  is isomorphic to the lexicographic product of the complete graph $K_{3^{3n-2}}$ and $K_3 \square K_3$.

\begin{theorem}\label{main2_0}
If  $H_{3n}$ is the graph defined in Construction \ref{construction-two},
then $H_{3n}^{2n}$ is isomorphic to the lexicographic product
$K_{3^{3n-2}} [K_{3}\square K_3]$ of the complete graph $K_{3^{3n-2}}$ and $K_{3}\square K_3$.
\end{theorem}

\begin{pf}
By Lemma~\ref{lemma_G}, $G_{3n}^{2n-1}$ is isomorphic to $K_{3\star 3^{3n-2}}$. Let
$P_1$, $P_2$, $\ldots$, $P_{3^{3n-2}}$ be the partite sets of $G_{3n}^{2n-1}$.
For each $i \in \{1,\ldots,3^{3n-2}\}$,
let $Q_i$ be a subset of vertices of $H_{3n}$ defined by
\[ Q_{i}=P_i \times V(K_3).\]
Note that $|Q_i|=9$, and  $\{Q_i \mid 1\le i\le 3^{3n-2}\}$ is a partition of $V(H_{3n})$.
%See Figure~\ref{fig2} for an illustration.

To show that  $H_{3n}^{2n}$ is isomorphic to the lexicographic product
$K_{3^{3n-2}} [K_{3}\square K_3]$, it is sufficient to show the following claims.
\medskip

\noindent {\bf Claim 1.} For any vertex $x \in Q_i$ and $y \in Q_j$ with $i \neq j$,
we have $d_{H_{3n}} (x, y) \leq 2n$.

\medskip

\noindent {\bf Claim 2.} For each $i\in \{ 1,\ldots,3^{3n-2}\}$, the subgraph of $H_{3n}^{2n}$ induced by $Q_i$ is
isomorphic to $K_3 \square K_3$.

\medskip

Since $V(H_{3n})=V(G_{3n})\times V(K_3)$, for any two vertices $x$ and $y$ in $H_{3n}$, we can denote
$x = (x_1, x_2)$ and $y = (y_1, y_2) $ with
$x_1,y_1\in V(G_{3n})$ and $x_2,y_2\in V(K_3)$.
By Lemma~\ref{lem3}, $d_{G_{3n}}(x_1,y_1)\le 2n$ for any $x_1, y_1 \in V(G_{3n})$ and $d_{G_{3n}}(x_1,y_1) = 2n$  only when $x_1 \neq y_1$ and $x_1 \sim y_1$.
It is well-known fact that
\[ d_{G_{3n}\square K_3} (x,y) = d_{G_{3n}}(x_1,y_1) +d_{K_3}(x_2,y_2).\]
Thus \[ d_{H_{3n}} (x,y)  \le 2n +1\]
and $d_{H_{3n}} (x,y)  = 2n+1$ if and only if
$d_{G_{3n}}(x_1,y_1)=2n$ and $d_{K_3}(x_2,y_2)=1$.

If $x\in Q_i$ and $y \in Q_j$ with $i \neq j$, then $d_{G_{3n}} (x_1, y_1 ) \leq 2n-1$ by Lemma \ref{lemma_G}.  Thus
\[ d_{H_{3n}} (x,y) \le 2n.\]
 This completes the proof of Claim 1.

Suppose that $x, y \in Q_i$ for some integer $i\in \{1,2,\ldots,3^{3n-2}\}$.
Then $x_1,y_1 \in P_i$ and $x_1 \sim y_1$.
Therefore $d_{H_{3n}} (x,y) = 2n+1$ if and only if $x_1 \neq  y_1$ and $x_2 \neq y_2$.
This means that two vertices $x$ and $y$ of $Q_i$ are  non-adjacent in $H_{3n}^{2n}$ if and only if
$x_1 \neq y_1$ and $x_2 \neq y_2$, which implies that
the subgraph of $H_{3n}^{2n}$ induced by $Q_i$ is isomorphic to $K_3 \square K_3$.
This completes the proof of Claim 2.

By Claim 1, for any vertex $x \in Q_i$ and $y \in Q_j$ with $i \neq j$, two vertices $x$ and $y$ are adjacent in $H_{3n}^{2n}$.
Therefore  $H_{3n}^{2n}$ is isomorphic to the lexicographic product
$K_{3^{3n-2}} [K_{3}\square K_3]$ by Claim 2.
\end{pf}

%%%%%%%%%%%%%%%%%%%%%%%%%%%%%%%%%%%%%%%%%%%%%%%%%%%%%%%%%%%
\begin{corollary}\label{main2}
For every integer $k \geq 2$ and for any positive integer $s$, there is a graph $G$ such that $G^k$ is isomorphic to the lexicographic product $K_{3^{3sk-2}} [K_{3}\square K_3]$  of the complete graph $K_{3^{3sk-2}}$ and $K_{3}\square K_3$.
\end{corollary}

\begin{pf}
For any positive integers $k$ and $s$, let $G = H_{3sk}^{2s}$.
Since
\[G^{k}=(H_{3sk}^{2s})^{k} =H_{3sk}^{2sk},\]
$G^{k}$ is isomorphic to
$K_{3^{3sk-2}}[K_{3}\square K_3]$ by Theorem~\ref{main2_0}.
\end{pf}

%%%%%%%%%%%%%%%%%%%%%%%%%%%%%%%%%%%%%%%%%%%%%%%%%%%%%%%%%%%
Next, we will compute the chromatic number and the list chromatic number of the lexicographic product $G[H]$ of $G$ and $H$
to complete the proof of Theorem \ref{main}.
For the chromatic number of the lexicographic product $G[H]$ of $G$ and $H$,
the following is known.

\begin{prop}  [\cite{GS1975}] \label{lexico}
If $\chi(H)=l$, then $\chi(G[H]) =\chi (G[K_l])$ for any graph $G$.
\end{prop}

\begin{lemma}\label{chromatic}
If $H = K_{3^{3n-2}} [K_{3}\square K_3]$ is the lexicographic product
of the complete graph $K_{3^{3n-2}}$ and $K_{3}\square K_3$, then
$\chi(H)=3^{3n-1}$.
\end{lemma}

\begin{pf}
Note that $\chi(K_3 \square K_3) =3$.
By  Proposition~\ref{lexico}, we have $\chi(H)= \chi(K_{3^{3n-2}}[K_3])$.
Since $K_{3^{3n-2}}[K_3]$ is a complete graph of order $3 ^{3n-1}$,
we have   $\chi(H)=3^{3n-1}$.
\end{pf}

\begin{lemma}\label{list_chromatic}
If $H = K_{3^{3n-2}} [K_{3}\square K_3]$ is the lexicographic product
of the complete graph $K_{3^{3n-2}}$ and $K_{3}\square K_3$, then
\[\chi_l(H) \ge \frac{10}{9}\cdot 3^{3n-1} -1.\]
\end{lemma}

\begin{pf}
Let $V(K_{3^{3n-2}}) = \{1, \ldots, 3^{3n-2}\}$, and denote
$V(K_{3^{3n-2}} [K_{3}\square K_3]) =  \{(x, y) : x \in \{1, \ldots, 3^{3n-2}\} \text { and } y \in  V(K_{3}\square K_3) \}$.

%Let $V(K_3) = \{0, 1, 2 \}$, and let
%$W_1 = \{(0, 0), (0, 1), (0, 2) \}$,
%$W_2 = \{(1, 0), (1, 1), (1, 2) \}$, and
%$W_3 = \{(2, 0), (2, 1), (2, 2) \}$.  Then $V(K_{3}\square K_3) = W_1 \cup W_2 \cup W_3$.  Note that for each %$1 \leq i \leq 3$, the subgraph of $K_{3}\square K_3$ induced by $W_i$ is a complete graph.

Let $V(K_3) = \{0, 1, 2 \}$, and we partition $V(K_{3}\square K_3)$ into three subsets $W_1, W_2,$ and $W_3$ such that
$W_1 = \{(0, 0), (0, 1), (0, 2) \}$,
$W_2 = \{(1, 0), (1, 1), (1, 2) \}$, and
$W_3 = \{(2, 0), (2, 1), (2, 2) \}$.
Note that for each $1 \leq i \leq 3$, the subgraph of $K_{3}\square K_3$ induced by $W_i$ is a complete graph.

For each $1 \leq i \leq 3$, let
\[ R_i= V(K_{3^{3n-2}}) \times W_i.\]
Then it is clear that $R_i$ is a subset of $V(H)$ and  $\{ R_1, R_2, R_3 \}$ is a partition of $V(H)$.
Note that the subgraph of $H$ induced by $R_i$ is a complete graph.

Let $t$ be an even integer.
Let $A_1$, $A_2$, $A_3$ be mutually disjoint sets such that $|A_1|=|A_2|=|A_3|=\frac{t}{2}$, and
let $A=A_1\cup A_2\cup A_3$.
For each vertex $v\in R_i$ with $1 \leq i \leq 3$, we define the list $L(v) = A\setminus A_i$.
We will show that $H$ is not $L$-choosable if $t < \frac{10}{9}\cdot 3^{3n-1}$.

First, we define a family of subsets of $V(H)$, denoted by $\{S_i : 1 \leq i \leq 3^{3n-2} \}$.  For each $i \in \{ 1, \ldots,  3^{3n-2}\}$, let
\[ S_i = \{(i, y): y \in V(K_{3}\square K_3) \}.\]
Then $\{S_i : 1 \leq i \leq 3^{3n-2} \}$ is a partition of $V(H)$.
Note that when $i$ and $j$ are distinct, for any two vertices $u \in S_i$ and $v \in S_j$,
$u$ and $v$ are adjacent in $H$ by the definition of the lexicographic product $K_{3^{3n-2}} [K_{3}\square K_3]$.

We will show that if there is a proper coloring $\phi$ such that $\phi(v)\in L(v)$ for each vertex $v$,
then $|\{ \phi(v) : v \in S_i \}|\ge  5$ for each $i  \in \{ 1, \ldots, 3^{3n-2} \}$.
Let $S_i$ be a subset of $H$ with $i \in \{1, \ldots, 3^{3n-2}\}$.
If there exist $u, v, w \in S_i$ such that
$\phi(u)=\phi(v)=\phi(w)$, then
$\{u,v,w\}$ forms an independent set of $H$.
Since each $R_j$ induces a complete graph in $H$,
any two of $u$, $v$, $w$ cannot belong to the same $R_j$.
Thus $L(u)\cap L(v)\cap L(w) = \emptyset$, which is a contradiction.
Hence, each color must be used at most two times in $S_i$.  Thus
$|\{ \phi(v) \mid v\in S_i \}|\ge \lceil\frac{9}{2}\rceil=5$.
Therefore $|\{ \phi(v) \mid v\in S_i \}|\ge  5$ for each $i  \in \{ 1, \ldots, 3^{3n-2} \}$.

Note that
\[A \supseteq \bigcup_{i=1}^{3^{3n-2}} \{ \phi(v) : v \in S_i\}.\]

Note that when $i$ and $j$ are distinct,
 for any two vertices $u \in S_i$ and $v \in S_j$,
$u$ and $v$ are adjacent.
Thus $\{ \phi(v) : v\in S_i\}$ and $\{ \phi(v) : v\in S_j \}$ are disjoint for any distinct $i$ and $j$.
%Let $u \in S_i$ and $v \in S_j$ be any vertices.  Then
%by the definition of the lexicographic product $K_{3^{3n-2}} [K_{3}\square K_3]$,
%$u$ and $v$ are adjacent in $H$ when $i$ and $j$ are distinct.
%Thus $\{ \phi(v) : v\in S_i\}$ and $\{ \phi(v) : v\in S_j \}$ are disjoint for any distinct $i$ and $j$.
Hence
\[  \frac{3t}{2} = |A|  \geq \big| \bigcup_{i=1}^{3^{3n-2}} \{ \phi(v) : v \in S_i\} \big|  \ge 5 \cdot3^{3n-2}.\]
Therefore
\[t \ge \frac{10}{3} \cdot 3^{3n-2} =\frac{10}{9} \cdot 3^{3n-1} .\]
This implies that for any even integer $t$ less than $\frac{10}{9} \cdot 3^{3n-1}$, we have $\chi_l(H)> t$. Since there is a possibility that $\chi_l(H)$ is odd, one can say that  $\chi_l(H) \ge \frac{10}{9} \cdot 3^{3n-1} - 1$.
\end{pf}

%%%%%%%%%%%%%%%%%%%%%%%%%%%%%%%%%%%%%%%%%%%%%%%%%%%%%%%%%%%%%%%%%%%%%%%%%%%%%%%%%%

%%%%%%%%%%%%%%%%%%%%%%%%%%%%%%%%%%%%%%%%%%%%%%%%%%%%%%%%%%%%%%%%

When $n = sk$, we have
$\chi_l (K_{3^{3sk-2}} [K_{3}\square K_3]) \geq \frac{10}{9}\cdot 3^{3sk-1} -1$ by Lemma~\ref{list_chromatic}
and $\chi(K_{3^{3sk-2}} [K_{3}\square K_3]) = 3^{3sk-1} $ by Lemma~\ref{chromatic}.
Therefoe Theorem~\ref{main} holds by
Corollary~\ref{main2} and Lemmas~\ref{chromatic} and~\ref{list_chromatic}.

%%%%%%%%%%%%%%%%%%%%%%%%%%%%%%%%%%%%%%%%%%%%%%%%%%%%%%%%%%%%%%%%%%%

%%%%%%%%%%%%%%%%%%%%%%%%%%%%%%%%%%%%%%%%%%%%%%%%%%%%%%%%%%%%%%%%%%%%

\section{Remark}

Since the List Square Coloring Conjecture is not true in general, a natural problem is to find the upper bound of the list chromatic number of $G^2$ for any graph $G$.
By observing some straightforward bounds on $\omega(G^2)$ and $\Delta(G^2)$, one immediately obtains $\chi_l(G^2) \le (\chi(G^2))^2$.
 Noel \cite{Noel} proposed the following two problems.

\begin{question} \label{question-one} [Noel \cite{Noel}]
Is there a function $f(x) = o(x^2)$ such that for every graph $G$, 
\[\chi_l (G^2) \leq f (\chi(G^2))?
\]
\end{question}

The example of Kim and Park shows that the function $f$ in Question \ref{question-one} must satisfy $f(x) = \Omega(x\log x)$. Noel \cite{Noel} asked whether it is possible
to obtain a general upper bound of the same order of magnitude.

\begin{question} \label{question-two} [Noel \cite{Noel}]
Does there exist a constant $c$ such that every graph satisfies
\[\chi_l (G^2) \leq c \chi(G^2) \log \chi(G^2)?
\]
\end{question}

In this paper, we show that there is no constant $k$ such that $G^k$ is chromatic-choosable for every graph $G$.  On the other hand, Gravier and Maffray \cite{GM98}  conjectured  that every claw-free graph is chromatic-choosable. As a relaxation of the conjecture, it is an interesting problem to answer the following question.

\begin{question} \label{claw-free}
Is there a constant $k$ such that $G^k$ is chromatic-choosable if $G$ is claw-free?
\end{question}

One could expect that if $G$ is a chromatic-choosable, then $G^2$ is also chromatic-choosable.  But, there exists a graph $G$ such that $G^2$ is not chromatic-choosable even though $G$ is chromatic-choosable.  The smallest example of Kim and Park's \cite{KP2013} is one of the graphs that have such property.
Therefore it is not clear whether $G^k$ is chromatic choosable for every $k \geq k_0$ even though $G^{k_0}$ is chromatic-choosable.
Thus by assumption that Gravier and Maffray's conjecture is true, it is an interesting problem to answer the following question.

\begin{question} \label{claw-free-k}
Is $G^k$  chromatic-choosable for every integer $k \geq 2$ if $G$ is claw-free?
\end{question}

%Thus giving an answer to Question \ref{claw-free} is an interesting problem even though

\medskip
%%%%%%%%%%%%%%%%%%%%%%%%%%%%%%%%%%%%%%%%%%%%%%%%%%%%%%%%%%%%%%%%%%%%%%
\noindent {\bf Acknowledgement.} The first and the second authors are supported by the National Research Foundation of Korea(NRF) grant funded by the Korea government (MEST)
 (No. 2011-0009729) and (No. 2010-0022142), respectively.
The third author was supported by the National Institute for Mathematical Sciences (NIMS) grant funded by the Korea government (B21203).

%%%%%%%%%%%%%%%%%%%%%%%%%%%%%%%%%%%%%%%%%%%%%%%%%%%%%%%%%%%%%%%%%%%%

 \end{document}